\documentstyle[11pt]{article}
\def\z{\bar z}
\def\o{\omega}
\def\d{\delta}
\date{ }
\begin{document}
\title{
 Willmore Submanifolds in a sphere
\thanks{The author is supported by a research fellowship of the Alexander von
Humboldt Stiftung 2001/2002 and the Zhongdian grant of the NSFC.}}
\author{Haizhong Li}
\maketitle

\begin{center}
\begin{minipage}{120mm}

\vskip 0.1in
\baselineskip 0.2 in
{{\bf Abstract:} Let $x:M\to S^{n+p}$ be an $n$-dimensional submanifold
in an $(n+p)$-dimensional unit sphere $S^{n+p}$, $x:M\to S^{n+p}$ is called a
Willmore submanifold if
it is a extremal submanifold to the following Willmore functional:
$$
\int_M(S-nH^2)^{\frac{n}{2}}dv,
$$
where $S=\sum\limits_{\alpha,i,j}(h^\alpha_{ij})^2$ is the square of the length
of the second fundamental form, $H$ is the mean curvature 
of $M$. In [13], author proved an integral inequality of Simons' type for 
$n$-dimensional compact Willmore hypersurfaces in $S^{n+1}$ and gave a characterization of {\it Willmore tori}.
In this paper, we generalize this result to $n$-dimensional compact 
Willmore submanifolds in $S^{n+p}$. In fact, we obtain an integral 
inequality of Simons' type for compact Willmore
submanifolds in $S^{n+p}$ and give a characterization of
{\it Willmore tori} and {\it Veronese surface} by use of our integral inequality.}
\medskip\noindent



\end{minipage}
\end{center}

\section *{1. Introduction}

Let $M$ be an $n$-dimensional compact submanifold of an $(n+p)$-dimensional
unit sphere space $S^{n+p}$. If $h^\alpha_{ij}$ denotes the second fundamental
form of $M$, $S$ denotes the square of the length of the second
fundamental form, $\bf{H}$ denotes the mean curvature vector and $H$
denotes the mean curvature of $M$, then
we have
$$
S=\sum\limits_{\alpha,i,j}(h^\alpha_{ij})^2,\quad
{\bf H}=\sum\limits_\alpha H^\alpha e_\alpha,\quad
H^\alpha=\frac{1}{n}\sum\limits_k h^\alpha_{kk},
\quad H=|{\bf H}|,
$$
where $e_{\alpha}$ ($n+1\leq\alpha\leq n+p$) are orthonormal normal
vector fields of $M$ in $S^{n+p}$.

       We define the following non-negative function on $M$
$$
\rho^2=S-nH^2,
\eqno (1.1)
$$
which vanishes exactly at the umbilic points of $M$.

Willmore functional is the following non-negative functional (see
[4], [21]
or [27])
$$
\int_M\rho^ndv=
\int_M (S-nH^2)^{\frac{n}{2}}dv,
$$
it was shown in [4], [21] and [27] that this functional is an invariant under
Moebius (or conformal) transformations of $S^{n+p}$. We use the term
{\it Willmore submanifolds} to call its critical points, because when $n=2$, 
the functional essentially coincides with the well-known Willmore functional 
$W(x)$
and its critical points are
the Willmore surfaces. Willmore conjecture says that $W(x)\geq 4\pi^2$ holds for all immersed tori $x:M\to S^3$. The conjecture was approached by Willmore [30], Li-Yau [18], Montiel-Ros [20], Ros [23], Langer-Singer [10] and many others(see [29], [31] and references there).

  In this paper, we first prove the following theorem (c.f.
Guo-Li-Wang [8], Pedit-Willmore [21]))

\medskip
{\bf Theorem 1.1} {\it Let $M$ be an $n$-dimensional submanifold in an
$(n+p)$-dimensional unit sphere $S^{n+p}$. Then $M$ is a Willmore
submanifold if and only if for $n+1\leq\alpha\leq n+p$
$$
\begin{array}{lcl}
&{}&
-\rho^{n-2}[SH^\alpha+\sum\limits_{\beta,i,j}H^\beta h^\beta_{ij}h^\alpha_{ij}
-\sum\limits_{\beta,i,j,k}h^\alpha_{ij}
h^\beta_{ik}h^\beta_{kj}-nH^2H^\alpha]
+(n-1)H^\alpha \Delta(\rho^{n-2})\\
&{}&\qquad
+2(n-1)\sum\limits_i (\rho^{n-2})_{i}H^\alpha_{,i}+(n-1)\rho^{n-2}
\Delta^\perp
H^\alpha
-\sum\limits_{i,j}(\rho^{n-2})_{i,j}(nH^\alpha\delta_{ij}
-h^\alpha_{ij})=0,
\end{array}
\eqno (1.2)
$$
where $\Delta (\rho^{n-2})=\sum\limits_i (\rho^{n-2})_{i,i}$,
$\Delta^\perp H^\alpha=\sum\limits_i H^\alpha_{,ii}$, and $(\rho^{n-2})_{i,j}$ 
is the Hessian of $\rho^{n-2}$ with respect to the induced metric $dx\cdot dx$, $H^\alpha_{,i}$  and $H^\alpha_{,ij}$ are the components of the first and second covariant derivative of the mean curvature vector field $\bf H$ ( see (2.14)-(2.17) ).}

\medskip
{\it Remark 1.1.} When $n=2$ and $p=1$, Theorem 1.1 was well-known (see 
Blaschke [1], Thomsen [26], Bryant [2] and chapter 7 of [31]).
When $n=2$ and $p\geq 1$ , Theorem 1.1 was proved by J. Weiner in [28], in this case
$(1.2)$ reduces to the following well-known equation of
Willmore surfaces (see [28] or [14])
$$
\Delta^\perp H^\alpha+\sum\limits_{\beta,i,j}h^\alpha_{ij}
h^\beta_{ij}H^\beta-2H^2H^\alpha=0,\quad
3\leq\alpha\leq 2+p.
\eqno (1.3)
$$
When $n\geq 2$ and $p=1$, Theorem 1.1 was proved by the author in [13].

\medskip
{\it Remark 1.2.} We should note that for $n\geq 2$, C.P. Wang [27] got the
Euler-Lagrange equation of Willmore functional for compact $n$-dimensional 
submanifolds without umbilical points in an $(n+p)$-dimensional unit sphere $S^{n+p}$
in terms of Moebius geometry. We also note that a different version of Theorem 1.1 was announced in Pedit-Willmore [21].
From the expression of $(1.2)$, in the cases $n=3$ and $n=5$ we need assume that $M$ has no umbilical points to guarantee $(\rho^{n-2})_{i,j}$ is continuous on $M$. We will make this assumption in this paper.

In order to state our main result, we first give the following
two important examples

{\bf Example 1} (see [13] or [8]). The tori
$$
W_{m,n-m}=S^m\left (\sqrt{\frac{n-m}{n}}\right )
\times S^{n-m}\left (\sqrt{\frac{m}{n}}\right ), \quad 1\leq m\leq n-1
\eqno (1.4)
$$
are Willmore hypersurfaces in $S^{n+1}$. We call $W_{m,n-m}$, $1\leq m\leq n-1$,
{\it Willmore tori}.
In fact, the principal curvatures $k_1,\cdots,k_n$ of $W_{m,n-m}$ are
$$
k_1=\cdots=k_m=\sqrt{\frac{m}{n-m}},\quad
k_{m+1}=\cdots=k_n=-\sqrt{\frac{n-m}{m}}.
\eqno (1.5)
$$
We have from $(1.5)$
$$
H={\frac{1}{n}}(m{\sqrt\frac{m}{n-m}}-(n-m){\sqrt\frac{n-m}{m}}),\quad
S=\frac{m^2}{n-m}+\frac{(n-m)^2}{m},
$$
$$
\sum\limits_{i,j,k}h_{ij}h_{jk}h_{ki}=\sum\limits_ik^3_i=
m(\frac{m}{n-m})^{\frac{3}{2}}-(n-m)(\frac{n-m}{m})^{\frac{3}{2}},
$$
where $h_{ij}=h^{n+1}_{ij}$. Thus we easily check that $(1.2)$  holds, i.e., $W_{m,n-m}$ are
Willmore hypersurfaces. In particular, we note that $\rho^2$ of $W_{m,n-m}$
for all $1\leq m\leq n-1$
satisfy 

$$\rho^2=n.
\eqno (1.6)
$$

We recall that well-known Clifford minimal
tori are 
$$
C_{m,n-m}=S^m\left (\sqrt{\frac{m}{n}}\right )\times S^{n-m}
\left (\sqrt{\frac{n-m}{n}}\right ),\qquad  1\leq m\leq n-1,
\eqno (1.7)
$$ 
It is remarkable that
a Willmore torus coincides with a Clifford minimal torus if and only if $n=2m$
for some $m$.

\medskip
{\it Remark 1.3.} When $n=2$, we can see from $(1.3)$ that all minimal surfaces
are Willmore surfaces. 
In [22], Pinkall constructed many compact non-minimal flat Willmore surfaces in $S^3$. In [3], Castro-Urbano constructed many compact non-minimal
Willmore surfaces in $R^4$. Ejiri [7] and Li-Vrancken [17] constructed many non-minimal flat tori in $S^5$ and $S^7$. Bryant [2] classified the Willmore spheres in $S^3$ and Montiel [19] classified the Willmore spheres in $S^4$.
When $n\geq 3$, minimal submanifolds are not
Willmore submanifolds in general, for example, Clifford minimal tori
$C_{m,n-m}=S^m\left ({\sqrt{\frac{m}{n}}}\right )\times S^{n-m}\left ({\sqrt\frac{n-m}{n}}\right )$ 
are not Willmore submanifolds when $n\not= 2m$. In [8], the authors proved that all $n$-dimensional minimal Einstein submanifolds in a sphere are Willmore submanifolds (we note that this result was stated in [21] by Pedit and Willmore).

\medskip
{\bf Example 2} (see [6] or [11]). {\it Veronese surface}. Let $(x,y,z)$ be the natural
coordinate system in $R^3$ and $u=(u_1,u_2,u_3,u_4,u_5)$ the natural
coordinate system in $R^5$. We consider the mapping defined by
$$
u_1=\frac{1}{\sqrt{3}}yz,\quad
u_2=\frac{1}{\sqrt{3}}xz,\quad
u_3=\frac{1}{\sqrt{3}}xy,\quad
$$
$$
u_4=\frac{1}{2\sqrt{3}}(x^2-y^2),\quad
u_5=\frac{1}{6}(x^2+y^2-2z^2),
$$
where $x^2+y^2+z^2=3$.
This defines an isometric immersion of $S^2(\sqrt{3})$ into $S^4(1)$.
Two points $(x,y,z)$ and $(-x,-y,-z)$ of $S^2(\sqrt{3})$ are mapped
into the same point of $S^4$. This real projective plane imbedded in
$S^4$ is called the {\it Veronese surface}.
We know that Veronese surface is a minimal
surface in $S^4$ (see [6] or [11]), thus it is a Willmore surface. We also
note that $\rho^2$ of the Veronese surface satisfies
$$
\rho^2=\frac{4}{3}.
\eqno (1.8)
$$

In the theory of minimal submanifolds in $S^{n+p}$,
the following J. Simons' integral inequality is well-known 

{\bf Theorem 1.2} (Simons [25], Lawson [11], Chern-Do Carmo-Kobayashi [6])
{\it Let
 $M$ be an $n$-dimensional ($n\geq 2$) compact minimal submanifold
in $(n+p)$-dimensional unit sphere $S^{n+p}$. Then we have
$$\int_M S\left (\frac{n}{2-1/p}-S\right )dv\leq 0.
\eqno (1.9)
$$
In particular, if
$$
0\leq S\leq \frac{n}{2-1/p},
\eqno (1.10)
$$
then either $S\equiv 0$ and $M$ is totally geodesic, or $S\equiv \frac{n}{2-1/p}$.
In the latter case, either $p=1$ and $M$ is a Clifford torus $C_{m,n-m}$,
or $n=2$, $p=2$ and $M$ is the Veronese surface.}

\medskip

In this paper we prove the following integral inequality of Simons' type for 
compact Willmore submanifolds in $S^{n+p}$. 

\medskip
{\bf Theorem 1.3} {\it Let $M$ be an $n$-dimensional ($n\geq 2$) compact Willmore 
submanifold in $(n+p)$-dimensional unit sphere $S^{n+p}$. Then we have
$$\int_M\rho^n\left (\frac{n}{2-1/p}-\rho^2\right )dv\leq 0.
\eqno (1.11)
$$
In particular, if
$$
0\leq \rho^2\leq \frac{n}{2-1/p},
\eqno (1.12)
$$
then either $\rho^2\equiv 0$ and $M$ is totally umbilic, or $\rho^2\equiv \frac{n}
{2-1/p}$. In the latter case, either $p=1$ and $M$ is a Willmore
torus $W_{m,n-m}$ defined by $(1.4)$; or $n=2$, $p=2$ and $M$
is the Veronese surface.}

\medskip
{\it Remark 1.4}. In case $p=1$, Theorem 1.3 was proved by the author
in [13]. We also classified all isoparametric Willmore hypersurfaces
in $S^{n+1}$ in [13].

\section *{2. Preliminaries}

Let $x:M\to S^{n+p}$ be an $n$-dimensional submanifold in an
$(n+p)$-dimensional unit sphere $S^{n+p}$. Let $\{e_1,\cdots,e_{n}\}$
be a local orthonormal
basis of $M$ with respect to the induced metric, $\{\theta_1,\cdots,
\theta_n\}$ are their dual form. Let $e_{n+1},\cdots, e_{n+p}$
be the local unit orthonormal normal
vector field. In this paper we make the following convention on the range of
indices:

$$
1\leq i,j,k\leq n; \qquad
n+1\leq \alpha,\beta,\gamma\leq n+p.
$$

Then we have the structure equations
$$
dx=\sum\limits_i \theta_i e_i,
\eqno (2.1)
$$
$$
de_i=\sum\limits_j\theta_{ij}e_j+\sum\limits_{\alpha,j} h^\alpha_{ij}
\theta_je_\alpha-\theta_i x,
\eqno (2.2)
$$
$$
de_\alpha=-\sum\limits_{i,j}h^\alpha_{ij}\theta_j e_i+\sum
\limits_\beta \theta_{\alpha\beta}e_\beta.
\eqno (2.3)
$$
The Gauss equations are
$$
R_{ijkl}=(\delta_{ik}\delta_{jl}-\delta_{il}\delta_{jk})
+\sum\limits_\alpha (h^\alpha_{ik}h^\alpha_{jl}-h^\alpha_{il}
h^\alpha_{jk}),
\eqno (2.4)
$$
$$
R_{ik}=(n-1)\delta_{ik}+n\sum\limits_\alpha H^\alpha h^\alpha_{ik}
-\sum\limits_{\alpha,j} h^\alpha_{ij}h^\alpha_{jk},
\eqno (2.5)
$$
$$
n(n-1)R=n(n-1)+n^2H^2-S,
\eqno (2.6)
$$
where $R$ is the normalized scalar curvature of $M$ and
$S=\sum_{\alpha,i,j}(h^\alpha_{ij})^2$ is the norm square of the
second fundamental form,
${\bf H}=\sum_\alpha H^\alpha e_\alpha=\frac{1}{n}\sum_\alpha
(\sum_k h^\alpha _{kk})e_\alpha$ is the mean curvarture vector
and $H=|{\bf H}|$ is the mean curvature of $M$.

The Codazzi equations are

$$
h^\alpha_{ijk}=h^\alpha_{ikj},
\eqno (2.7)
$$
where the covariant derivative of $h^\alpha_{ij}$ is defined by 
$$
\sum\limits_k h^\alpha_{ijk}\theta_k=dh^\alpha_{ij}
+\sum\limits_kh^\alpha_{kj}\theta_{ki}
+\sum\limits_kh^\alpha_{ik}\theta_{kj}+\sum\limits_\beta
h^\beta_{ij}\theta_{\beta\alpha}.
\eqno (2.8)
$$

The second covariant derivative of $h^\alpha_{ij}$ is defined by 
$$
\sum\limits_l h^\alpha_{ijkl}\theta_l=dh^\alpha_{ijk}
+\sum\limits_l h^\alpha_{ljk }\theta_{li}
+\sum\limits_l h^\alpha_{ilk }\theta_{lj}
+\sum\limits_l h^\alpha_{ijl }\theta_{lk}
+\sum\limits_\beta h^\beta_{ijk}\theta_{\beta\alpha}.
\eqno (2.9)
$$

By exterior differentiation of $(2.8)$, we have the following
Ricci identities
$$
h^\alpha_{ijkl}-h^\alpha_{ijlk}=\sum\limits_m h^\alpha_{mj}R_{mikl}
+\sum\limits_m h^\alpha_{im}R_{mjkl}+\sum\limits_\beta h^\beta_{ij}
R_{\beta\alpha kl}.
\eqno (2.10)
$$

The Ricci equations are
$$
R_{\alpha\beta ij}=\sum\limits_k (h^\alpha_{ik}h^\beta_{kj}-
h^\beta_{ik}h^\alpha_{kj}).
\eqno (2.11)
$$

We define the following non-negative function on $M$
$$
\rho^2=S-nH^2,
\eqno (2.12)
$$
which vanishes exactly at the umbilical points of $M$.

    Willmore functional is the following functional (see [4], [21] and [27])
$$
W(x):=\int_M\rho^ndv=\int_M (S-nH^2)^{\frac{n}{2}}dv.
\eqno (2.13)
$$
By Gauss equation $(2.6)$, $(2.13)$ is equivalent to
$$
W(x)=[n(n-1)]^{\frac{n}{2}}\int_M (H^2-R+1)^{\frac{n}{2}}dv.
\eqno (2.13)'
$$
It was shown
in [4], [21] and [27] that this functional is an invariant under
Moebius (or conformal) transformations of $S^{n+p}$. We use the term
{\it Willmore submanifolds} to call its critical points. When $n=2$, the functional essentially
coincides with the well-known Willmore functional and its critical points are
{\it Willmore surfaces}.

We define the first, second covariant derivatives and Laplacian of the mean curvature vector field ${\bf H}=\sum\limits_\alpha H^\alpha e_\alpha$ in the normal bundle $N(M)$ as follows

$$
\sum\limits_i H^\alpha_{,i}\theta_i=dH^\alpha+\sum\limits_\beta H^\beta\theta_{\beta\alpha},
\eqno (2.14)
$$
$$
\sum\limits_j H^\alpha_{,ij}\theta_j=dH^\alpha_{,i}+
\sum\limits_j H^\alpha_{,j}\theta_{ji}+\sum\limits_\beta H^\beta_{,i}\theta_{\beta\alpha},
\eqno (2.15)
$$

$$
\Delta^\perp H^\alpha=\sum\limits_i H^\alpha_{,ii},\qquad H^\alpha=\frac{1}{n}
\sum\limits_k h^\alpha_{kk}.
\eqno (2.16)
$$

Let $f$ be a smooth function on $M$, we define the first, second covariant 
derivatives $f_i$, $f_{i,j}$ and Laplacian of $f$ as follows
$$
df=\sum\limits_i f_i\theta_i,\qquad
\sum\limits_j f_{i,j}\theta_j=df_i+\sum\limits_jf_j\theta_{ji},\qquad \Delta f=
\sum\limits_i f_{i,i}.
\eqno (2.17)
$$

\section *{3.  Proof of Theorem 1.1}
\vskip 5mm\noindent

 Let $x:M\to S^{n+p}$ be a compact submanifold in a unit sphere
$S^{n+p}$.  From Wang [27], we have the following relations of the connections of
the Mobeius metric ${\tilde\rho}^2dx\cdot dx$ and induced metric $dx\cdot dx$
$$
\omega_{ij}=\theta_{ij}+({\rm ln}{\tilde\rho})_i\theta_j-({\rm ln}{\tilde\rho})_j\theta_i,
\quad
\omega_{\alpha\beta}=\theta_{\alpha\beta},\quad
{\tilde\rho}=\sqrt{\frac{n}{n-1}}\rho,
\eqno (3.1)
$$
where $\rho$ is defined by $(1.1)$.

We get by use of $(3.15)$ of [27] and $(3.1)$
$$
\begin{array}{lcl}
\sum\limits_j{\tilde\rho} C^\alpha_{i,j}\theta_j&=& 
\sum\limits_j C^\alpha_{i,j}\omega_j\\
&=&dC^\alpha_i+\sum\limits_j C^\alpha_j\omega_{ji}+
\sum\limits_\beta C^\beta_i\omega_{\beta\alpha}\\
&=&dC^\alpha_i+\sum\limits_j C^\alpha_j\theta_{ji}+
\sum\limits_\beta C^\beta_i\theta_{\beta\alpha}+
\sum\limits_k [({\rm ln}{\tilde\rho})_k\theta_i-({\rm ln}{\tilde\rho})_i
\theta_k]\cdot C^\alpha_k,
\end{array}
$$
thus
$$
\begin{array}{lcl}
{\tilde\rho} C^\alpha_{i,j}&=&
2{\tilde\rho}^{-3}{\tilde\rho}_j[H^\alpha_{,i}+\sum\limits_k
(h^\alpha_{ik}-H^\alpha\delta_{ik})
({\rm ln}{\tilde\rho})_k]-{\tilde\rho}^{-2}[H^\alpha_{,ij}+\sum\limits_k
(h^\alpha_{ikj}
-H^\alpha_{,j}\delta_{ik})({\rm ln}{\tilde\rho})_k]\\
&{}&\quad
-{\tilde\rho}^{-2}\sum\limits_k (h^\alpha_{ik}-H^\alpha\delta_{ik})({\rm ln}{\tilde\rho})_{k,j}+
\sum\limits_k C^\alpha_k({\rm ln}{\tilde\rho})_k\delta_{ij}-({\rm ln}{\tilde\rho})_iC^\alpha_j,
\end{array}
\eqno (3.2)
$$
where $\{({\rm ln}{\tilde\rho})_{k,j}\}$ is the Hessian-Matrix of 
${\rm ln}{\tilde\rho}$ with respect to the induced metric $dx\cdot dx$ and
$\{H^\alpha_{,i}\}$ is defined by $(2.14)$.

Letting $i=j$, making summation over $i$ in $(3.2)$ and using $(3.15)$ of [27], we have
$$
\begin{array}{lcl}
\sum\limits_i C^\alpha_{i,i}
&=&
-{\tilde\rho}^{-3}\Delta^\perp H^\alpha -{\tilde\rho}^{-3}
\sum\limits_{i,k}(h^\alpha_{ik}-H^\alpha\delta_{ik})
({\rm ln}{\tilde\rho})_{k,i}\\
&{}&\quad
-2(n-2){\tilde\rho}^{-3}\sum\limits_i H^\alpha_{,i}({\rm ln}{\tilde\rho})_i
-(n-3){\tilde\rho}^{-3}\sum\limits_{i,j}({\rm ln}
{\tilde\rho})_i({\rm ln}{\tilde\rho})_j(h^\alpha_{ij}-H^\alpha\delta_{ij}),
\end{array}
\eqno (3.3)
$$
where $\Delta^\perp H^\alpha$ is defined by $(2.16)$.

On the other hand, we have from $(3.10)$ and $(3.14)$ of [27]
$$
\begin{array}{lcl}
&{}&
\sum\limits_{i,j}A_{ij}B^\alpha_{ij}
+\sum\limits_{\beta,i,j,k}B^\beta_{ik}B^\beta_{kj}B^\alpha_{ij}\\
&=&
-{\tilde\rho}^{-3}\sum\limits_{i,j}({\rm ln}{\tilde\rho})_{i,j}(h^\alpha_{ij}
-H^\alpha\delta_{ij})
+{\tilde\rho}^{-3}\sum\limits_{i,j}({\rm ln}{\tilde\rho})_i({\rm ln}{\tilde\rho})_j
(h^\alpha
_{ij}-H^\alpha \delta_{ij})\\
&{}&
\quad
+{\tilde\rho}^{-3}[\sum\limits_{\beta,i,j,k}h^\beta_{ik}h^\beta_{kj}h^\alpha_{ij}-\sum\limits_{\beta,i,j}H^\beta h^\beta_{ij}
h^\alpha_{ij}-S H^\alpha+ nH^2H^\alpha].
\end{array}
\eqno (3.4)
$$

Putting $(3.3)$ and $(3.4)$ into the following Willmore condition (see
$(2.34)$ and $(4.27)$ of [27])
$$
-(n-1)\sum\limits_iC^\alpha_{i,i}
+\sum\limits_{i,j}A_{ij}B^\alpha_{ij}
+\sum\limits_{\beta,i,j,k}
B^\alpha_{ij} B^\beta_{ik}B^\beta_{kj}=0,
\eqno (3.5)
$$
we get
$$
\begin{array}{lcl}
&{}&
\frac{n-1}{\rho^{n+1}}\{-\frac{\rho^{n-2}}{n-1}
[SH^\alpha+\sum\limits_{\beta,i,j}H^\beta h^\beta_{ij}h^\alpha_{ij}
-\sum\limits_{\beta,i,j,k}h^\alpha_{ij}
h^\beta_{ik}h^\beta_{kj}-n H^2H^\alpha]\\
&{}&\quad
+\rho^{n-2}\Delta^\perp H^\alpha+\frac{n-2}{n-1}\rho^{n-2}
\sum\limits_{i,j}({\rm ln}\rho)_{i,j}(h^\alpha_{ij}-H^\alpha\delta_{ij})
+2(n-2)
\rho^{n-2}\sum\limits_i({\rm ln}\rho)_i H^\alpha_{,i}\\
&{}&\quad
+\frac{(n-2)^2}{n-1}\rho^{n-2}
\sum\limits_{i,j}({\rm ln}\rho)_i({\rm ln}\rho)_j(h^\alpha_{ij}-H^\alpha \delta_{ij})\}
=0.
\end{array}
\eqno (3.5)'
$$

We can check the following identity by a direct computation
$$
\begin{array}{lcl}
&{}&-\frac{1}{n-1}\sum\limits_{i,j}(\rho^{n-2})_{i,j}(nH^\alpha\delta_{ij}-h^\alpha_{ij})
+\rho^{n-2}\Delta^\perp H^\alpha +2\sum\limits_i(\rho^{n-2})_iH^\alpha_{,i}
+H^\alpha\Delta(\rho^{n-2})\\
&=&\frac{(n-2)^2}{n-1}\rho^{n-2}\sum\limits_{i,j}({\rm ln}\rho)_i({\rm ln}\rho)_j
(h^\alpha_{ij}-H^\alpha\delta_{ij})+
\frac{n-2}{n-1}\sum\limits_{i,j}\rho^{n-2}({\rm ln}\rho)_{i,j}(h^\alpha_{ij}-H^\alpha
\delta_{ij})\\
&{}&\quad
+2(n-2)\rho^{n-2}\sum\limits_i({\rm ln}\rho)_i H^\alpha_{,i}+\rho^{n-2}\Delta^\perp
H^\alpha.
\end{array}
\eqno (3.6)
$$
Thus $(3.5)'$ is equivalent to $(1.2)$ by use of $(3.6)$. We
complete the proof of Theorem 1.1.

\medskip\noindent

{\it Remark 3.1} Fix index $\alpha$ with $n+1\leq \alpha\leq n+p$,
define $\Box^\alpha:M\to R$ by
$$
\Box^\alpha f=(nH^\alpha\delta_{ij}-h^\alpha_{ij})f_{i,j},
$$
where $f$ is any smooth function on $M$ and $f_{i,j}$ is defined by $(2.17)$. We know that $\Box^\alpha$
is a self-adjoint operator (cf. Cheng-Yau [5], Li [15,16]).
It is remarkable that this operator naturally
appears in Willmore equation $(1.2)$. In fact,
by use of this self-adjoint operator, Willmore equation $(1.2)$ can be 
written as the following equivalent form

$$
\begin{array}{lcl}
&{}&
-\rho^{n-2}[SH^\alpha+\sum\limits_{\beta,i,j}H^\beta h^\beta_{ij}h^\alpha_{ij}
-\sum\limits_{\beta,i,j,k}h^\alpha_{ij}
h^\beta_{ik}h^\beta_{kj}-nH^2H^\alpha] +(n-1)\rho^{n-2}\Delta^\perp H^\alpha\\
&{}&\qquad
+2(n-1)\sum\limits_i (\rho^{n-2})_iH^\alpha_{,i}+(n-1)H^\alpha \Delta(\rho^{n-2})-\Box^\alpha(\rho^{n-2})=0,\quad n+1\leq\alpha\leq n+p.
\end{array}
\eqno (1.2)'
$$

\section *{4. The Lemmas }

We first prove the following Lemma (c.f. Simons [25],
Chern-Do Carmo-Kobayashi [6] or Schoen-Simon-Yau [24])

\medskip
{\bf Lemma 4.1} {\it Let $M$ be an $n$-dimensional ($n\geq 2$)
submanifold in $S^{n+p}$. Then we have
$$
\begin{array}{lcl}
\frac{1}{2}\Delta\rho^2&=& |\nabla h|^2-n^2|\nabla^\perp {\bf H}|^2
+\sum\limits_{\alpha,i,j,k} (h^\alpha_{ij}h^\alpha_{kki})_j
+n\sum\limits_{\alpha,\beta,i,j,m} H^\beta h^\beta_{mj}h^\alpha_{ij}
h^\alpha_{im}\\
&{}&\quad
+n\rho^2
-\sum\limits_{\alpha,\beta,i,j,k,m} h^\alpha_{ij}h^\beta_{ij}h^\alpha_{mk}h^\beta_{mk}
-\sum\limits_{\alpha,\beta,j,k}
(R_{\beta\alpha jk})^2-\frac{1}{2}\Delta (nH^2),
\end{array}
\eqno (4.1)
$$
where $|\nabla h|^2=\sum\limits_{\alpha,i,j,k} (h^\alpha_{ijk})^2$ and  $|\nabla^\perp
{ \bf H}|^2=\sum\limits_{\alpha,i}(H^\alpha_{,i})^2$, $H^\alpha_{,i}$ is defined by $(2.14)$. }

\medskip
{\it Proof}. By the definition of $\Delta $ and $\rho^2$, we have by
use of $(2.7)$ and $(2.10)$
$$
\begin{array}{lcl}
\frac{1}{2}\Delta\rho^2&=&\frac{1}{2}\Delta(\sum\limits_{\alpha,i,j}
 (h^\alpha_{ij})^2)
-\frac{1}{2} \Delta (n H^2)\\
&=&
\sum\limits_{\alpha,i,j,k} (h^\alpha_{ijk})^2+\sum\limits_{\alpha,i,j,k}
 h^\alpha_{ij}h^\alpha_{kijk}-
\frac{1}{2}\Delta (n H^2)\\
&=&|\nabla h|^2-n^2|\nabla ^\perp{\bf H}|^2+
\sum\limits_{\alpha,i,j,k} (h^\alpha_{ij}h^\alpha_{kki})_j
+\sum\limits_{\alpha,i,j,k,m}
h^\alpha_{ij} h^\alpha_{mk}R_{mijk}\\
&{}&\quad+
\sum\limits_{\alpha,i,j,m} h^\alpha_{ij}h^\alpha_{im}R_{mj}
+\sum\limits_{\alpha,\beta,i,j,k} h^\alpha_{ij}
h^\beta_{ik} R_{\beta\alpha jk}-\frac{1}{2}\Delta (n H^2).
\end{array}
\eqno (4.2)
$$

By use of $(2.4)$ and $(2.5)$, we have
$$
\begin{array}{lcl}
&{}&
\sum\limits_{\alpha,i,j,k,m} h^\alpha_{ij}h^\alpha_{mk}R_{mijk}
+\sum\limits_{\alpha,i,j,m} h^\alpha_{ij}
h^\alpha_{im}R_{mj}
+\sum\limits_{\alpha,\beta,i,j,k} h^\alpha_{ji}h^\beta_{ik}R_{\beta\alpha jk}\\
&=& nS-n^2H^2-\sum\limits_{\alpha,\beta,i,j,k,m} h^\alpha_{ij}h^\beta_{ij}h^\alpha_{mk}h^\beta_{mk}
+n\sum\limits_{\alpha,\beta,i,j,m} H^\beta h^\beta_{mj}h^\alpha_{ij}h^\alpha_{im}\\
&{}&\quad
-[\sum\limits_{\alpha,\beta,i,j,m,l} h^\alpha_{ij}h^\alpha_{im}h^\beta_{ml}h^\beta_{lj}
-\sum\limits_{\alpha,\beta,i,j,k,m} h^\alpha_{ij}h^\alpha_{km}h^\beta_{jm}h^\beta_{ik}-\sum\limits_{\alpha,\beta,i,j,k}
h^\alpha_{ji}h^\beta_{ik}R_{\beta\alpha jk}].
\end{array}
\eqno (4.3)
$$

On the other hand, we have by $(2.11)$
$$
\begin{array}{lcl}
\sum\limits_{\alpha,\beta,j,k} (R_{\beta\alpha jk})^2&=&
\sum\limits_{\alpha,\beta,i,j,k} (h^\beta_{ji}h^\alpha_{ik}-h^\alpha_{ji}h^\beta_{ik})
R_{\beta\alpha jk}\\
&=&
\sum\limits_{\alpha,\beta,i,j,k,l} h^\beta_{ji}h^\alpha_{ik}(h^\beta_{jl}h^\alpha_{lk}-
h^\beta_{kl}h^\alpha_{lj})
-\sum\limits_{\alpha,\beta,i,j,k} h^\alpha_{ji}h^\beta_{ik}
R_{\beta\alpha jk}\\
&=&
\sum\limits_{\alpha,\beta,i,j,m,l} h^\alpha_{ij}h^\alpha_{im}h^\beta_{ml}h^\beta_{lj}
-\sum\limits_{\alpha,\beta,i,j,k,m} h^\alpha_{ij}h^\alpha_{km}h^\beta_{jm}h^\beta_{ik}-\sum\limits_{\alpha,\beta,i,j,k}
h^\alpha_{ji}h^\beta_{ik}R_{\beta\alpha jk}.
\end{array}
\eqno (4.4)
$$

Putting $(4.3)$ into $(4.2)$, we obtain $(4.1)$ by use of $(4.4)$.

\medskip
{\bf Lemma 4.2} {\it Let $M$ be an $n$-dimensional ($n\geq 2$)
submanifold in $S^{n+p}$, then we have
$$
|\nabla h|^2\geq \frac{3n^2}{n+2}|\nabla^\perp {\bf H}|^2,
\eqno (4.5)
$$
where $|\nabla h|^2=\sum\limits_{\alpha,i,j,k}(h^\alpha_{ijk})^2$,
$|\nabla^\perp {\bf H}|^2=\sum\limits_{\alpha,i}(H^\alpha_{,i})^2$,
$H^\alpha_{,i}$ is defined by $(2.14)$.}

{\it Proof}.  We construct the
following symmetric trace-free tensor (c.f. [9] or [13])
$$
F^\alpha_{ijk}
=
h^\alpha_{ijk}-
\frac{n}{n+2}(H^\alpha_i\delta_{jk}+H^\alpha_j\delta_{ik}
+H^\alpha_k\delta_{ij}).
\eqno (4.6)
$$
Then we can easily compute that
$$
|F|^2=\sum\limits_{\alpha,i,j,k}(F^\alpha_{ijk})^2
=|\nabla h|^2-\frac{3n^2}{n+2}|\nabla^\perp {\bf H}|^2.
$$
Then we have
$$
|\nabla h|^2\geq \frac{3n^2}{n+2}|\nabla^\perp {\bf H}|^2,
$$
which proves the Lemma 4.2.

Define tensors 
$$
{\tilde h}^\alpha_{ij}=h^\alpha_{ij}-H^\alpha\delta_{ij},
\eqno (4.7)
$$
$$
{\tilde\sigma}_{\alpha\beta}=\sum_{i,j}{\tilde h}^\alpha_{ij}
{\tilde h}^\beta_{ij},\quad
{\sigma}_{\alpha\beta}=\sum_{i,j}{ h}^\alpha_{ij}
{ h}^\beta_{ij}.
\eqno (4.8)
$$

Then the $(p\times p)$-matrix  $({\tilde\sigma}_{\alpha\beta})$ is
symmetric and can be assumed to be diagonized for a suitable choice of
$e_{n+1},\cdots, e_{n+p}$. We set
$$
{\tilde\sigma}_{\alpha\beta}={\tilde\sigma}_\alpha\delta_{\alpha\beta}.
\eqno (4.9)
$$

We have by a direct calculation
$$
\sum\limits_k {\tilde h}^\alpha_{kk}=0,\quad
{\tilde \sigma}_{\alpha\beta}=\sigma_{\alpha\beta}-nH^\alpha H^\beta,
\quad
\rho^2=\sum\limits_\alpha{ \tilde \sigma}_{\alpha}= S-nH^2,
\eqno (4.10)
$$

$$
\sum\limits_{\beta,i,j,k}h^\alpha_{ij}h^\beta_{ik}h^\beta_{kj}
=\sum\limits_{\beta,i,j,k}{\tilde h}^\alpha_{ij}
{\tilde h}^\beta_{ik}{\tilde h}^\beta_{kj}+2
\sum\limits_{\beta,i,j}H^\beta{\tilde h}^\beta_{ij}
{\tilde h}^\alpha_{ij} +H^\alpha\rho^2+nH^\alpha H^2,
\eqno (4.11)
$$

$$
\sum\limits_{\alpha,i,j,m}h^\beta_{mj}h^\alpha_{ij}h^\alpha_{im}
=\sum\limits_{\alpha,i,j,m}{\tilde h}^\beta_{mj}
{\tilde h}^\alpha_{ij}{\tilde h}^\alpha_{im}+2
\sum\limits_{\alpha,i,j}H^\alpha{\tilde h}^\alpha_{ij}
{\tilde h}^\beta_{ij} +H^\beta\rho^2+nH^2 H^\beta.
\eqno (4.12)
$$

From $(4.7),(4.10),(4.11)$ and Theorem 1.1, we have

{\bf Lemma 4.3}. {\it Let $M$ be an $n$-dimensional
submanifold in an $(n+p)$-dimensional unit sphere $S^{n+p}$.
Then $M$ is a Willmore submanifold if and only if for
$n+1\leq\alpha\leq n+p$

$$
\begin{array}{lcl}
&{}&
(n-1)\rho^{n-2}\Delta^\perp H^\alpha \\
&=&
-2(n-1)\sum\limits_i (\rho^{n-2})_i H^\alpha_{,i}-(n-1)H^\alpha\Delta(\rho^{n-2})\\
&{}&\quad
-\rho^{n-2}(\sum\limits_\beta H^\beta{\tilde\sigma}_{\alpha\beta}
+\sum\limits_{\beta,i,j,k}
{\tilde h}^\alpha_{ij}{\tilde h}^\beta_{ik}{\tilde h}^\beta_{kj})
+\sum\limits_{i,j}(\rho^{n-2})_{i,j}
(nH^\alpha\delta_{ij}-h^\alpha_{ij}),
\end{array}
\eqno (4.13)
$$
where $(\rho^{n-2})_{i,j}$, $\Delta (\rho^{n-2})$ and 
$\Delta^\perp H^\alpha$ are defined by $(2.17)$ and $(2.16)$.} 

In general, for a matrix $A=(a_{ij})$ we denote by $N(A)$ the square of
the norm of $A$, i.e.,
$$
N(A)={\rm trace} (A\cdot A^t)=\sum\limits_{i,j}(a_{ij})^2.
$$
Clearly, $N(A)=N(T^{t}AT)$ for any orthogonal matrix $T$.

\medskip
We need the following lemma (see Chern-Do Carmo-Kobayashi [6])

\medskip
{\bf Lemma 4.4}. {\it Let $A$ and $B$ be symmetric $(n\times n)$-matrices.
Then
$$
N(AB-BA)\leq 2N(A)\cdot N(B),
$$
and the equality holds for nonzero matrices $A$ and $B$ if and only if
$A$ and $B$ can be transformed simultaneously by an orthogonal matrix into
multiples of ${\tilde A}$ and ${\tilde B}$ respectively, where
$$
{\tilde A}=
\left (
\begin{array}{ccccc}
0&1&0&\cdots&0\\
1&0&0&\cdots&0\\
0&0&0&\cdots&0\\
\vdots&\vdots&\vdots&\ddots&\vdots\\
0&0&0&\cdots&0
\end{array}
\right )
\qquad
{\tilde B}=
\left (
\begin{array}{ccccc}
1&0&0&\cdots&0\\
0&-1&0&\cdots&0\\
0&0&0&\cdots&0\\
\vdots&\vdots&\vdots&\ddots&\vdots\\
0&0&0&\cdots&0
\end{array}
\right ).
$$
Moreover, if $A_1, A_2$ and $A_3$ are $(n\times n)$-symmetric matrices and if
$$
N(A_\alpha A_\beta -A_\beta A_\alpha)=2N(A_\alpha)\cdot N(A_\beta),\quad
1\leq \alpha,\beta\leq 3,
$$
then at least one of the matrices $A_\alpha$ must be zero.}

\medskip
{\bf Lemma 4.5}. {\it Let $x:M\to S^{n+p}$ be an $n$-dimensional submanifold
in $S^{n+p}$. Then}
$$
\begin{array}{lcl}
\frac{1}{2}\Delta\rho^2&\geq & |\nabla h|^2-n^2|\nabla^\perp {\bf H}|^2
+\sum\limits_{\alpha,i,j,k}(h^\alpha_{ij}h^\alpha_{kki})_j
+n\sum\limits_{\alpha,\beta,i,j,m} H^\beta {\tilde h}^\beta_{mj}{\tilde h}^\alpha_{ij}{\tilde h}^
\alpha_{im}\\
&{}&{\displaystyle
\quad
+n\rho^2
+nH^2\rho^2-(2-\frac{1}{p})\rho^4
-\frac{1}{2}\Delta (nH^2)}.
\end{array}
\eqno (4.14)
$$

{\it Proof}.  By use of $(4.10)$ and $(4.12)$, $(4.1)$ becomes
$$
\begin{array}{lcl}
\frac{1}{2}\Delta\rho^2&=& |\nabla h|^2-n^2|\nabla^\perp {\bf H}|^2
+\sum\limits_{\alpha,i,j,k}(h^\alpha_{ij}h^\alpha_{kki})_j
+n\sum\limits_{\alpha,\beta,i,j,m} H^\beta {\tilde h}^\beta_{mj}{\tilde h}^\alpha_{ij}{\tilde h}^\alpha_{im}\\
&{}&\quad
+n\rho^2
+nH^2\rho^2
-\sum\limits_{\alpha,\beta}{\tilde \sigma}^2_{\alpha\beta}
-\sum\limits_{\alpha,\beta,j,k}
(R_{\beta\alpha jk})^2-\frac{1}{2}\Delta (nH^2).
\end{array}
\eqno (4.15)
$$

By $(2.11)$, we have 
$$
\begin{array}{lcl}
\sum\limits_{\alpha,\beta,j,k} (R_{\beta\alpha jk})^2
&=&\sum\limits_{\alpha,\beta,j,k}(\sum\limits_l h^\beta_{jl}h^\alpha_{lk}
-\sum\limits_l h^\alpha_{jl}h^\beta_{lk})^2\\
&=&\sum\limits_{\alpha,\beta,j,k}(\sum\limits_l {\tilde h}^\beta_{jl}
{\tilde h}^\alpha_{lk}
-\sum\limits_l {\tilde h}^\alpha_{jl}{\tilde h}^\beta_{lk})^2\\
&=& \sum\limits_{\alpha,\beta}N(A_\alpha A_\beta-A_\beta A_\alpha),
\end{array}
\eqno (4.16)
$$
where
$$
A_\alpha:=({\tilde h}^\alpha_{ij})
=(h^\alpha_{ij}-H^\alpha\delta_{ij}).
\eqno (4.17)
$$

By use of Lemma 4.4, $(4.9)$, $(4.10)$ and $(4.16)$, we get
$$
\begin{array}{lcl}
-\sum\limits_{\alpha,\beta}{\tilde\sigma}_{\alpha\beta}^2-
\sum\limits_{\alpha,\beta,j,k}(R
_{\beta\alpha jk})^2
&=&
-\sum\limits_{\alpha}{\tilde \sigma}^2_\alpha-\sum\limits_
{\alpha,\beta}N(A_\alpha A_\beta-A_\beta A_\alpha)\\
&\geq&
-\sum\limits_{\alpha}{\tilde \sigma}^2_\alpha
-2\sum\limits_{\alpha\not=\beta}{\tilde \sigma}_\alpha
{\tilde\sigma}_\beta\\
&=&
-2(\sum\limits_\alpha{\tilde\sigma}_\alpha)^2+\sum\limits_\alpha
{\tilde\sigma}_\alpha^2\\
&\geq&{\displaystyle
-2\rho^4+\frac{1}{p}
(\sum\limits_\alpha{\tilde\sigma}_\alpha)^2 }\\
&=&{\displaystyle
-(2-\frac{1}{p})\rho^4. }
\end{array}
\eqno (4.18)
$$

Putting $(4.18)$ into $(4.15)$, we get $(4.14)$.

\medskip
{\bf Lemma 4.6}. {\it Let $M$ be an $n$-dimensional submanifold
in an $(n+p)$-dimensional unit sphere $S^{n+p}$, then we have}

$$
\begin{array}{lcl}
\frac{1}{2}\Delta(\rho^n)
&\geq&
\frac{1}{2} n\{-n(n-1)
\rho^{n-2}
|\nabla^\perp{\bf H}|^2
+\rho^{n-2}
\sum\limits_{\alpha,i,j,k}(h^\alpha_{ij}h^\alpha_{kki})_j
+\rho^n[n+nH^2\\
&{}&\quad
-(2-\frac{1}{p})\rho^2] 
+n\rho^{n-2}
\sum\limits_{\alpha,\beta,i,j,m} H^\beta{\tilde h}^\beta_{mj}{\tilde h}
^\alpha_{ji}{\tilde h}^\alpha
_{im} -\frac{1}{2}\rho^{n-2}\Delta(nH^2)\}.
\end{array}
\eqno (4.19)
$$

{\it Proof.} First it is easy to check the following calculation

$$
\begin{array}{lcl}
\frac{1}{2}\Delta(\rho^n)
&=&\frac{1}{2}\Delta[(\rho^2)^{\frac{n}{2}}]
=\frac{1}{2}n(n-2)\rho^{n-2}\sum\limits_i \rho_i^2+
\frac{n}{4}\rho^{n-2}\Delta(\rho^2)\\
&\geq& \frac{n}{4}\rho^{n-2}\Delta(\rho^2).
\end{array}
\eqno (4.20)
$$

Noting
$$
\begin{array}{lcl}
|\nabla h|^2-n^2|\nabla^\perp {\bf H}|^2
&=&
(|\nabla h|^2-\frac{3n^2}{n+2}|\nabla^\perp {\bf H}|^2)
+(\frac{3n^2}{n+2}-n)|\nabla^\perp {\bf H}|^2-n(n-1)|\nabla^\perp
{\bf H}|^2\\
&\geq&
(|\nabla h|^2-\frac{3n^2}{n+2}|\nabla^\perp {\bf H}|^2)
-n(n-1)|\nabla^\perp
{\bf H}|^2,
\end{array}
\eqno (4.21)
$$
we get from  $(4.5)$, $(4.14)$ and $(4.21)$
$$
\begin{array}{lcl}
\frac{1}{2}\Delta \rho^2&\geq & -n(n-1)|\nabla^\perp {\bf H}|^2
+\sum\limits_{\alpha,i,j,k} (h^\alpha_{ij}h^\alpha_{kki})_j
+n\sum\limits_{\alpha,\beta,i,j,m} H ^\beta {\tilde h}^\beta_{mj}{\tilde h}^\alpha_{ij}{\tilde h}^
\alpha_{im}\\
&{}&{\displaystyle
\quad
+n\rho^2
+nH^2\rho^2-(2-\frac{1}{p})\rho^4
-\frac{1}{2}\Delta (nH^2)}.
\end{array}
\eqno (4.22)
$$

We obtain $(4.19)$ by putting $(4.22)$ into $(4.20)$.

\medskip
{\bf Lemma 4.7}. {\it Let $x:M\to S^{n+p}$ be an $n$-dimensional
compact submanifold in $S^{n+p}$. Let $f$ and $g$ be
two any smooth functions on $M$.} Then we have
$$
\int_M g \Delta fdv=-\int_M(\sum\limits_if_ig_i)dv=\int_Mf\Delta g dv.
\eqno (4.23)
$$

{\it Proof}. Integrating the following identities over $M$
$$
g\Delta f=g\sum\limits_i f_{i,i}=
\sum\limits_i (f_ig)_i-\sum\limits_i f_ig_i,
$$
$$
f\Delta g=f\sum\limits_i g_{i,i}=
\sum\limits_i (g_if)_i-\sum\limits_i f_ig_i,
$$
we get $(4.23)$.

\medskip
{\bf Lemma 4.8}. {\it Let $M$ be an $n$-dimensional compact Willmore 
submanifold in an $(n+p)$-dimensional unit sphere $S^{n+p}$,
then we have}

$$
\begin{array}{lcl}
&{}&-n(n-1)\int_M\rho^{n-2}|\nabla^\perp {\bf H}|^2dv
+n\int_M \rho^{n-2}(\sum\limits_{\alpha,\beta,i,j,m}H^\beta
{\tilde h}^\beta_{mj}{\tilde h}^\alpha_{ji}{\tilde h}^\alpha_{im})dv\\
&=&
-\frac{1}{2}n(n+1)\int_M\sum\limits_i (\rho^{n-2})_i(H^2)_idv
-n\int_M \rho^{n-2}\sum\limits_{\alpha,\beta} H^\alpha H^\beta
{\tilde \sigma}_{\alpha\beta}dv\\
&{}&\quad
-n\int_M (\sum\limits_{\alpha,i,j} H^\alpha h^\alpha_{ij}
(\rho^{n-2})_{i,j})dv.
\end{array}
\eqno (4.24)
$$
{\it Proof.} We have the following calculation by use of Stokes formula 
and Lemma 4.7
$$
\begin{array}{lcl}

&{}&-n(n-1)\int_M\rho^{n-2}|\nabla^\perp {\bf H}|^2dv
+n\int_M \rho^{n-2}(\sum\limits_{\alpha,\beta,i,j,m}H^\beta
{\tilde h}^\beta_{mj}{\tilde h}^\alpha_{ji}{\tilde h}^\alpha_{im})dv\\
&=&
\quad
\frac{1}{2}n(n-1)\int_M (\sum\limits_i (\rho^{n-2})_i(H^2)_i)dv
+n(n-1)\int_M\rho^{n-2}\sum\limits_\alpha H^\alpha \Delta^\perp H^\alpha dv\\
&{}&\quad
+
n\int_M \rho^{n-2}(\sum\limits_{\alpha,\beta,i,j,m}H^\beta
{\tilde h}^\beta_{mj}{\tilde h}^\alpha_{ji}{\tilde h}^\alpha_{im})dv.
\end{array}
\eqno (4.25)
$$

Putting $(4.13)$ into $(4.25)$, using Stokes formula and Lemma 4.7, 
we get $(4.24)$.

\medskip
{\bf Lemma 4.9}. {\it Let $M$ be an $n$-dimensional  
compact submanifold in an $(n+p)$-dimensional unit sphere $S^{n+p}$,
then we have}

$$
\int_M \rho^{n-2}\sum\limits_{\alpha,i,j,k}(h^\alpha_{ij}h^\alpha_
{kki})_jdv
=n\int_M(\sum\limits_{\alpha,i,j}
H^\alpha h^\alpha_{ij}(\rho^{n-2})_{i,j})dv
+\frac{n^2}{2}\int_M \sum\limits_i (\rho^{n-2})_i(H^2)_i dv.
\eqno (4.26)
$$

{\it Proof.}  We have the following calculation by use of Stokes formula and Lemma 4.7

$$
\begin{array}{lcl}
\int_M\rho^{n-2}\sum\limits_{\alpha,i,j,k}( h^\alpha_{ij}
h^\alpha_{kki})_jdv
&=&-\int_M\sum\limits_{\alpha,i,j,k}(\rho^{n-2})_jh^\alpha_{ij}
h^\alpha_{kki}dv \\
&=&-\int_M\sum\limits_{\alpha,i,j,k}((\rho^{n-2})_{j}h^\alpha_{ij}
h^\alpha_{kk})_idv+\int_M\sum\limits_{\alpha,i,j,k}
(\rho^{n-2})_{j,i}h^\alpha_{ij}h^\alpha_{kk}dv\\
&{}&\quad+n^2\int_M\sum\limits_{\alpha,j} H^\alpha(\rho^{n-2})_{j}
H^\alpha_{,j}dv\\
&=& n\int_M\sum\limits_{\alpha,i,j} H^\alpha h^\alpha_{ij}
(\rho^{n-2})_{i,j}dv
+\frac{1}{2}n^2 \int_M\sum\limits_i 
 (\rho^{n-2})_i(H^2)_i dv.
\end{array}
$$

By use of Stokes formula and Lemma 4.7, we also have 

\medskip
{\bf Lemma 4.10}. {\it Let $M$ be an $n$-dimensional  
compact submanifold in an $(n+p)$-dimensional unit sphere $S^{n+p}$,
then we have}

$$
-\frac{1}{2}\int_M \rho^{n-2}\Delta (nH^2)dv
=\frac{n}{2}\int_M\sum\limits_i (\rho^{n-2})_i(H^2)_idv.
\eqno (4.27)
$$

\section *{5. The Proof of Theorem 1.3 }

Now we begin to prove the Theorem 1.3.
 Integrating $(4.19)$ over $M$ and using Stokes formula, we have
$$
\begin{array}{lcl}
0&\geq&  
\frac{n}{2}\{-n(n-1)\int_M \rho^{n-2}|\nabla^\perp {\bf H}|^2dv
+n\int_M
\rho^{n-2}
\sum\limits_{\alpha,\beta,i,j,m}H^\beta{\tilde h}^\beta_{mj}
{\tilde h}^\alpha_{ji}{\tilde h}^\alpha_{im}dv\\
&{}&\quad
+\int_M \rho^{n-2}\sum\limits_{\alpha,i,j,k}
(h^\alpha_{ij}h^\alpha_{kki})_jdv -\frac{1}{2}\int_M\rho^{n-2}\Delta(nH^2)dv\\
&{}&\quad
+\int_M \rho^n(n+nH^2-(2-\frac{1}{p})\rho^2)dv\}.
\end{array}
\eqno (5.1)
$$

Putting $(4.24)$, $(4.26)$ and  $(4.27)$ into $(5.1)$, we get 

$$
\begin{array}{lcl}
0&\geq&{\displaystyle
\frac{n}{2}\{
-n\int_M \rho^{n-2}\sum\limits_{\alpha,\beta} H^\alpha H^\beta
{\tilde \sigma}_{\alpha\beta}dv
+\int_M \rho^n(n+nH^2-(2-\frac{1}{p})\rho^2)dv\} }\\
&=&
{\displaystyle
\frac{n}{2}\{
n\int_M \rho^{n-2}(H^2\rho^2-\sum\limits_{\alpha,\beta}
H^\alpha H^\beta {\tilde \sigma}_{\alpha\beta})dv
+
\int_M \rho^n(n-(2-\frac{1}{p})\rho^2)dv\}}\\
&\geq&
{\displaystyle
\frac{n}{2}\int_M \rho^n(n-(2-\frac{1}{p})\rho^2)dv },
\end{array}
\eqno (5.2)
$$
where we used
$$
\sum\limits_{\alpha,\beta}H^\alpha H^\beta
{\tilde\sigma}_{\alpha\beta}=\sum\limits_\alpha (H^\alpha)^2{\tilde
\sigma}_\alpha\leq \sum\limits_\alpha(H^\alpha)^2\cdot
\sum\limits_\beta
{\tilde\sigma}_\beta=H^2\rho^2.
\eqno (5.3)
$$

Thus we reach the following integral inequality of Simons' type
$$
\int_M\rho^n(n-(2-\frac{1}{p})\rho^2)dv\leq 0.
\eqno (5.4)
$$

Therefore we have proved the integral inequality $(1.11)$ in
Theorem 1.3.

If $(1.12)$ holds, then we conclude from $(5.4)$ that
either $\rho^2\equiv 0$, or $\rho^2\equiv n/(2-\frac{1}{p})$.
In the first case, we know that $S\equiv nH^2$, i.e.
$M$ is totally umbilic; in the latter case, i.e.,

$$
\rho^2=\sum\limits_{\alpha,i,j}({\tilde h}^\alpha_{ij})^2
\equiv n/(2-\frac{1}{p}),
\eqno (5.5)
$$
we have from $(4.21)$ (in this case, $(4.21)$ becomes an equality)
$$
(\frac{3n^2}{n+2}-n)|\nabla ^\perp{\bf H}|^2=0,
$$
we have
$\nabla ^\perp{\bf H}=0$, thus we have from $(4.5)$ (in this case
$(4.5)$ becomes an equality)
$$
\nabla h=0, \quad i.e., \quad
h^\alpha_{ijk}=0.
\eqno (5.6)
$$
It follows that all inequalities in $(4.18),(4.19)$, $(4.21)$, $(4.22)$,
$(5.1)$ and $(5.2)$ are actually equalities. In deriving $(4.19)$ from
$(4.18)$, we made use of inequalities
$N(A_\alpha A_\beta-A_\beta
A_\alpha)\leq2 N(A_\alpha)\cdot N(A_\beta)$. Hence

$$
N(A_\alpha A_\beta-A_\beta
A_\alpha)=2 N(A_\alpha)\cdot N(A_\beta),
\quad
\alpha\not=\beta.
\eqno (5.7)
$$

From $(4.18)$ we obtain
$$
{\tilde\sigma}_1={\tilde\sigma}_2=\cdots
={\tilde\sigma}_p,
\eqno (5.8)
$$
where ${\tilde\sigma}_\alpha$ is defined by $(4.9)$.

From $(5.7)$ and Lemma 4.4, we conclude that at most two of the
matrices $A_\alpha=({\tilde h}^\alpha_{ij})=(h^\alpha_{ij}-H^\alpha
\delta_{ij})$ are nonzero, in which case they can be assumed to be
scalar multiples of ${\tilde A}$ and ${\tilde B}$ in Lemma 4.4. We now
consider the case $p=1$ and $p\geq 2$ separately.

\medskip
{\bf Case $p=1$}. From $(5.6)$, we know that $M$ is an isoparametric
hypersurface.
In this case Theorem 1.3 already has been proved by the author in [13] (see
Theorem 3 of [13]). We conclude that  $M$ is one of Willmore tori,
that is, $M=W_{m,n-m}$ for some $m$ with $1\leq m\leq n-1$.

\medskip
{\bf Case $p\geq 2$}. In this case,
we know that at most two of $A_\alpha=({\tilde h}^\alpha_{ij})$,
$\alpha=n+1,\cdots,n+p$, are different from zero. If all of 
$A_\alpha=({\tilde h}^\alpha_{ij})$ are zero, which is contradiction with
$(5.5)$.
Assume that only
one of them, say $A_\alpha$, is different zero, which is contradiction with 
$(5.8)$.  Therefore we may 
assume that
$$
A_{n+1}=\lambda{\tilde A},\qquad
A_{n+2}=\mu{\tilde  B},\quad \lambda,\mu\not=0,
$$
$$
A_\alpha=0, \quad  \alpha\geq n+3,
$$
where ${\tilde A}$ and ${\tilde B}$ are defined in Lemma 4.4.

In this case, the inequality $(5.3)$ is actually equality, that is,
$$
\sum\limits_{\alpha,\beta}H^\alpha H^\beta
{\tilde\sigma}_{\alpha\beta}
=H^2\rho^2.
\eqno (5.9)
$$
In fact, $(5.9)$ is
$$
2\lambda^2(H^{n+1})^2+2\mu^2(H^{n+2})^2=
[(H^{n+1})^2+(H^{n+2})^2+\cdots+(H^{n+p})^2][2\lambda^2+2\mu^2].
$$
In view of $\lambda,\mu\not=0$, we follow
$$
H^\alpha=0,\quad n+1\leq\alpha\leq n+p,
$$
that is, ${\bf H}=0$, i.e., $M$ is a minimal submanfold in $S^{n+p}$,
$(5.5)$ becomes
$$
S=\sum\limits_{\alpha,i,j}(h^\alpha_{ij})^2=\frac{n}{2-1/p}.
$$
From Main Theorem of Chern-Do Carmo-Kobayashi [6], we conclude
that $n=2,p=2$ and $M$ is a Veronese
surface defined by example 2.
We complete the proof of Theorem 1.3.

\section *{6. Some Related Results}

For $p\geq 2$, we can improve Theorem 1.3 as follows

\medskip 
{\bf Theorem 6.1} {\it Let $M$ be an $n$-dimensional ($n\geq 2$) compact Willmore 
submanifold in $(n+p)$-dimensional unit sphere $S^{n+p}$. Then we have
$$\int_M\rho^n\left (\frac{2}{3}n-\rho^2\right )dv\leq 0.
\eqno (6.1)
$$
In particular, if
$$
0\leq \rho^2\leq \frac{2}{3}n,
\eqno (6.2)
$$
then either $\rho^2\equiv 0$ and $M$ is totally umbilic, or $\rho^2\equiv
\frac{2}
{3}n$. In the latter case, $n=2$, $p=2$ and $M$
is the Veronese surface.}

\medskip
{\it Remark 6.1} When $n=2$, Theorem 6.1 was proved by author in [14] (see Theorem 3 of [14]). When $p\geq 2$, the pinching constant $2n/3$ in Theorem 6.1 is better than
the pinching constant $n/(2-1/p)$ in Theorem 1.3.

\medskip
We need the following lemma to prove our Theorem 6.1

\medskip
{\bf Lemma 6.1} (see Theorem 1 of [12]). {\it Let $A_{n+1},\cdots, A_{n+p}$ be 
symmetric $(n\times n)$-matrices, which are defined by $(4.17)$. Set 
${\tilde\sigma}_{\alpha\beta}={\rm tr}(A^t_\alpha A_\beta)$,
${\tilde\sigma}_\alpha={\tilde\sigma}_{\alpha\alpha}=N(A_\alpha):=
{\rm tr}(A^t_\alpha A_\alpha)$, $\rho^2=\sum\limits_\alpha
{\tilde\sigma}_\alpha$ (see (4.8)-(4.10)), we have

$$
\sum\limits_{\alpha,\beta} N(A_\alpha A_\beta-A_\beta A_\alpha)
+\sum\limits_{\alpha,\beta}
{\tilde\sigma}_{\alpha\beta}^2\leq \frac{3}{2}\rho^4.
\eqno (6.3)
$$
Moreover, the equality holds if and only if at most two matrices $A_\alpha$
and $A_\beta$ are not zero and these two matrices can be transformed 
simultaneously by an orthogonal matrix into scalar multiples of ${\tilde A}$,
${\tilde B}$, respectively, where ${\tilde A}$ and ${\tilde B}$ are defined in
Lemma 4.4.}

\medskip
{\bf Proof of Theorem 6.1} In the proof of Lemma 4.5, using Lemma 6.1 instead
of Lemma 4.4, we can get

$$
\begin{array}{lcl}
\frac{1}{2}\Delta\rho^2&\geq & |\nabla h|^2-n^2|\nabla^\perp {\bf H}|^2
+\sum\limits_{\alpha,i,j,k}(h^\alpha_{ij}h^\alpha_{kki})_j
+n\sum\limits_{\alpha,\beta,i,j,m} H^\beta {\tilde h}^\beta_{mj}{\tilde h}^\alpha_{ij}{\tilde h}^
\alpha_{im}\\
&{}&{\displaystyle
\quad
+n\rho^2
+nH^2\rho^2-\frac{3}{2}\rho^4
-\frac{1}{2}\Delta (nH^2)}.
\end{array}
\eqno (6.4))
$$

Repeating process of the proof of Lemma 4.6, we have
$$
\begin{array}{lcl}
\frac{1}{2}\Delta(\rho^n)
&\geq&
\frac{1}{2} n\{-n(n-1)
\rho^{n-2}
|\nabla^\perp{\bf H}|^2
+\rho^{n-2}
\sum\limits_{\alpha,i,j,k}(h^\alpha_{ij}h^\alpha_{kki})_j
+\rho^n[n+nH^2\\
&{}&\quad
-\frac{3}{2}\rho^2] 
+n\rho^{n-2}
\sum\limits_{\alpha,\beta,i,j,m} H^\beta{\tilde h}^\beta_{mj}{\tilde h}
^\alpha_{ji}{\tilde h}^\alpha
_{im} -\frac{1}{2}\rho^{n-2}\Delta(nH^2)\}.
\end{array}
\eqno (6.5)
$$

Integrating $(6.5)$ over $M$, using Stokes formula, $(4.24), (4.26)$ and $(4.27)$
we can get $(6.1)$ by the similar argument as the proof of Theorem 1.3.

If $(6.2)$ holds, then we conclude from $(6.1)$ that either $\rho^2\equiv 0$
and $M$ is totally umbilic, 
or $\rho^2\equiv \frac{2}{3}n$. In the latter case, i.e.,
$$
\rho^2\equiv \frac{2}{3}n,
\eqno (6.6)
$$
using the similar argument as the proof of theorem 1.3, we can conclude that $\nabla^\perp {\bf H}\equiv 0$, $h^\alpha_{ijk}\equiv 0$ and ${\bf H}\equiv 0$, i.e., $M$ is a minimal submanifold with $(6.6)$. From Theorem 3 of [12], we know that $n=2$, $p=2$ and $M$ is the Veronese surface.
We complete the proof of Theorem 6.1.

\medskip
We now mention the following example

\medskip
{\bf Example 6.1} ([8]). $M_{m_1,\cdots,m_{p+1}}:=S^{m_1}(a_1)\times \cdots\times
S^{m_{p+1}}(a_{p+1})$ is an $n$-dimensional Willmore submanifold in $S^{n+p}$, where $n=m_1+\cdots+m_{p+1}$ and $a_i$ are defined by
$$
a_i=\sqrt{\frac{n-m_i}{np}},\quad i=1,\cdots,p+1.
$$
From Proposition 4.1 of [8], we can check that the function $\rho^2$ of
$M_{m_1,\cdots,m_{p+1}}$ satisfies
$\rho^2=np$ and 
$M_{m_1,\cdots,m_{p+1}}$ is a minimal submanifold if and only if 
$$m_1=m_2=\cdots=m_{p+1}=\frac{n}{p+1},\qquad
a_i=\sqrt{\frac{1}{p+1}},
\quad i=1,\cdots,p+1.
$$

Motivated by Chern's conjecture about minimal submanifolds in $S^{n+p}$ (see
[6]) and Theorem 1.3, we propose the following problem:

\medskip
{\bf Problem}. Let $x:M\to S^{n+p}$ be an $n$-dimensional compact Willmore
submanifold in an $(n+p)$-dimensional unit sphere $S^{n+p}$ with $\rho^2:=S-nH^2=$ constant. Is the set of values of $\rho^2$ discrete?

\vskip 0.2in

{\bf Acknowledgements:} The author began to do this research work
during his stay in department of mathematics of Harvard University
as a visitor in the academic year of 1999-2000.
He would like to express his thanks to Prof. S.-T.Yau for his
encouragement and help. He would like to express his thanks to
Prof. C.L.Terng for her useful suggestions and help.
He also thank the referee for some helpful comments which make 
this paper more readable than the original version of the paper.

\medskip\noindent
\begin{tabbing}
XXXXXXXXXXXXXXXXXXXXXXXXXX*\=\kill
{Permanent Address}\>{Present Address}\\
{Department of Mathematical Sciences}\>{Institut f\"ur Mathematik, MA 8-3}\\
{Tsinghua University}\>{Technische Universit\"at Berlin}\\
{100084, Beijing}\>Strasse des 17. Juni 136, 10623 Berlin\\
{People's Republic of China}\> Germany\\
{E-mail: hli@math.tsinghua.edu.cn}\>
{E-mail: hli@math.tu-berlin.de}
\end{tabbing}

\end {document}